\documentclass[12pt]{article}
\usepackage[T2A]{fontenc}
\usepackage[russian]{babel}
\usepackage{amsfonts}
\usepackage{amsmath}
\usepackage{amssymb}
\usepackage{amsthm}

\title{Сильно эллиптические операторы с сингулярными коэффициентами}

\author{
М.~И.~Нейман--заде\footnote{Работа поддержана грантом
РФФИ N 01-01-00691}, А.~А.~Шкаликов\footnote{Работа поддержана грантом
РФФИ N 00-15-96100}}

\newcommand\pa{\partial}

\newcommand\ov{\overline}

\newcommand\wt{\widetilde}
\newcommand\opn{\operatorname}

\renewcommand\Re{\opn{Re}}
\renewcommand\Im{\opn{Im}}
\renewcommand{\tg}{\opn{tg}}

\newcommand\D{{\mathcal D}}
\renewcommand\H{{\mathcal H}}

\newcommand\R{{\mathbb R}}

\newcommand\al{\alpha}
\newcommand\be{\beta}
\newcommand\de{\delta}

\newcommand\ga{\gamma}

\newcommand\De{\Delta}

\newcommand\la{\lambda}

\newcommand\si{\sigma}

\renewcommand\th{\theta}

\renewcommand\phi{\varphi}

\newcommand\ep{\varepsilon}

\newcommand{\urto}{\stackrel{R}{\Longrightarrow}}

\def\H{{\mathfrak H}}

\def\intl{\int\limits}
\def\suml{\sum\limits}

\newlength{\MMM}
\newlength{\MMMM}
\newcommand{\Mo}{%
\makebox{\parbox[l][\MMMM][b]{\MMM}{%
\raisebox{0.12em}{$\stackrel{\circ}{M}$}}}}

\newtheorem{theor}{Теорема}[section]
\newtheorem{theorem}[theor]{Теорема}

\newtheorem{lemm}[theor]{Лемма}

\newtheorem{remark}{Замечание}[section]
\newenvironment{proof*nodot}{\par\smallskip{\bf Доказательство}\rm}
			{\hfill $\square$\medskip}

\newenvironment{theorem*}[1]{\par\smallskip{\bf #1.}\it}{\par\medskip}

\begin{document}

\makeatletter


\renewcommand{\subsection}{\@startsection{subsection}{2}%
{0pt}{2.25ex plus 1ex minus .2ex}%
{0.5ex plus .2ex}{\large\bf}
}

\makeatother

\renewcommand{\bibname}{Список литературы}

\maketitle

В этой статье мы будем изучать сильно эллиптические
операторы
\begin{align}\label{III:eq:L}
&L=\suml_{|\al|,|\be|\le m} D^\al c_{\al,\be}(x) D^\be  \\
\notag
&c_{\al,\be}(x)\in \left\{
\begin{aligned}
L_\infty, &\quad|\al|=|\be|=m,\\
D',&\quad|\al|+|\be|<2m,
\end{aligned}
\right.
\end{align}
в $\R^n$ порядка $2m$. Используются стандартные обозначения:
$D^\al=i^{|\al|}\dfrac{\pa^{|\al|}}{\pa^{\al_1}x_1\dots
\pa^{\al_n}x_n}$, $\al=(\al_1,\dots,\al_n)$ и $\be$ --- мультииндексы,
$|\al|=\suml_{j=1}^n|\al_j|$.

Пусть также
$$
  L_0=\suml_{|\al|,|\be|=m} D^\al c_{\al,\be} D^\be
$$
--- главная часть оператора $L$.
Напомним (см., например,~\cite{LM},\cite{Fi}), что оператор $L_0$ называется
равномерно сильно эллиптическим, если
$$
 \Re\suml_{|\al|=|\be|=m} c_{\al,be}(x)\xi^\al \xi^\be\ge
\de|\xi|^{2m},\quad x\in\R^n,\quad \xi\in\R^n,
$$
где постоянная $\de$ не зависит от $x$, $\xi\in\R^n$.

Согласно теореме Иосиды (см.~\cite{Ios}), верен следующий факт:
{\it
если оператор $L_0$ равномерно сильно эллиптичен, а коэффициенты $c_{\al,\be}$
непрерывные функции, то выполняется коэрцитивная оценка
\begin{equation}\label{III:eq:Re>}
\Re (L_0 u,u)\ge\de \|\nabla^m u\|^2_{L_2(\R^n)}
\end{equation}
для всех $u\in S$, где $S$ --- пространство Шварца, а постоянная
$\de_1$ не зависит от $u$.
}
Эту оценку называют таккже неравенством
Гординга~\cite{Fi}.

В дальнейшем мы нигде не будем использовать непрерывность
коэффициентов оператора $L_0$, а только их принадлежность пространству
$L_\infty$ и оценку~\eqref{III:eq:Re>}. Поэтому в случае
$c_{\al,\be}\in L_\infty$
далее постулируется справедливость оценки~\eqref{III:eq:Re>}.

Из оценки~\eqref{III:eq:Re>}
и того, что коэффициенты
главной части $c_{\al,\be}$ существенно ограничены, следует, что
полуторалинейная форма
$$
l_0(u,v)=\intl_{\R^n}\suml_{|\al|=|\be|=m} c_{\al,\be} D^\be u(x)
\ov{D^\al v(x)}dx
$$
секториальна, замыкаема, и область определения ее замыкания
равна $H^m_2(\R^n)$. Поэтому, согласно первой теореме о представлении,
существует единственный $m$--секториальный оператор
в пространстве $L_2(\R^n)$, который мы обозначим также через $L_0$,
такой, что
$$
  l_0(u,v)=(L_0 u,v),\quad u\in \D(L_0),\quad v\in H^m_2(\R^n).
$$

В определении~\eqref{III:eq:L} мы предположили, что коэффициенты
$c_{\al,\be}(x)$ при $|\al|+|\be|<2m$ являются распределениями.
Основной вопрос, исследуемый в этой главе, --- найти
достаточные условия на коэффициенты $c_{\al,\be}$, гарантиующие корректное
определение
$L$ как секториального оператора в $L_2(\R^n)$.
Ответ на этот вопрос будет дан в терминах принадлежности коэффициентов
пространствам мультипликаторов из соболевского пространства $H_2^k$ в
$H_2^{-l}$. Поскольку точное описание таких пространств сложно, в финальной
теореме мы получим достаточные условия принадлежности обобщенных
функций этим пространствам в терминах их принадлежности соболевским
пространствам с негативными индексами гладкости.

Кроме того, мы показываем "физическую корректность"
оператора $L$. Именно,
если $c_{\al,\be,n}$, $|\al|+|\be|<m$ --- гладкие функции,
сходящиеся к $c_{\al,\be}$
по норме соответствующего пространства мультипликаторов, то
последовательность операторов $L_n$ будет сходиться к $L$ в смысле
равномерной резольвентной сходимости.

Нам будет удобнее сначала получить результаты в общей форме.

\subsection{Обобщенная сумма секториальных операторов. Приближения
операторных сумм}

Напомним, что оператор $T$ в гильбертовом пространстве $\H$
называется секториальным, если найдутся числа $\th<\pi/2$ и $M<\infty$
такие, что
$$
 |Im(Tx,x)|\le \tg\th(Re(Tx,x)+M(x,x)),\quad \forall x\in\H,\ \|x\|=1.
$$
Секториальный оператор $T$ называется $m$-секториальным, если
внешность сектора $|\Im\la|\le \tg\th(\Re\la+M)$ содержится в резольвентном
множестве оператора $T$.

Далее в этом параграфе мы будем предполагать, что $T$ --- секториальный
оператор в гильбертовом пространстве $\H$. Пусть $\H_1$ --- гильбертово
пространство, плотно вложенное в $\H$, а оператор $T$ таков, что
\begin{equation}\label{III:eq:ReTeq1}
c_1\|x\|_1^2\le\Re(Tx,x)\le c_2\|x\|_1^2\quad\forall x\in\D(T),
\end{equation}
где $\|\cdot\|_1$ --- норма в $\H_1$, а постоянные $c_1$ и $c_2$
не зависят от $x$.
В этом случае квадратичная форма $\Re(Tx,x)$ допускает замыкание
на все пространство $\H_1$, и, согласно первой теореме о предствалении
(см.~\cite[гл. VI.2]{Ka}), существует равномерно положительный самосопряженый
оператор $T_0$ такой, что
$$
(T_0x,x)=\Re(Tx,x)\qquad\text{при $x\in\D(T_0)\cap D(T)$},
$$
кроме того, $\D(T_0)\subset\H_1$ и $\D(T_0)$ плотно в $\H^1$. Известно также, 
что в этом случае $\D(T_0^{1/2})=\H_1$. Через $\H_{-1}$ обозначим дуальное
пространство к $\H_1$ по отношению к скалярному произведению в $\H$.
В частности, запись $(f,g)$ обозначает как скалярное произведение в $\H$, 
так и действие функционала $f\in\H_{-1}$ на элемент $g\in\H_1$.

\begin{theorem}\label{III:th:gensum}
Пусть $H_1$, $H$ --- пара гильбертовых пространств, причем
$\H_1$ плотно вложено в $\H$. Пусть $T$ --- секториальный оператор в $\H$, 
причем выполнена двусторонняя оценка~\eqref{III:eq:ReTeq1}. 
Пусть $Q$ --- ограниченный оператор из $\H_1$ в $\H_{-1}$, 
причем выполнена оценка
\begin{equation}\label{III:eq:Qsub}
|(Qf,f)|\le\ep(f,f)_1+M(\ep)(f,f)~\qquad\text{для некоторого $\ep<1$.}
\end{equation}
Тогда $T+Q$ --- замкнутый ограниченный оператор из пространства $\H_1$
в $\H_{-1}$, для которого форма $((T+Q)x,x)$, определенная на $\H_1$,
замкнуа и секториальная. Сужение оператора $T+Q$ на $\H$ (обозначим
его $S=T\wt+Q$ и назовем обобщенной суммой) есть $m$--секториальный оператор
в $\H$. Если последовательность операторов $Q_n:\H_1\to\H_{-1}$ такова, что
\begin{equation}\label{III:eq:Qconv}
\|Q_n-Q\|\to 0\qquad\text{при $n\to\infty$}
\end{equation}
(здесь берется операторная норма из $\H_1$ в $\H_{-1}$), то
\begin{equation}\label{III:eq:Sconv}
S_n=T\wt+Q_n\urto T\wt+Q=S.
\end{equation}
При этом спектры $\si(S_n)$ сходятся к спектру $\si(S)$ сверху.
\end{theorem}
\begin{proof}
Так как $T$ секториальный и выполнена оценка~\eqref{III:eq:ReTeq1},
то
$$
 c_1\|x_1\|^2_1\le|\Im(Tx,x)|\le c_2\|x\|_1^2
$$
с некоторыми константами $c_1$, $c_2>0$. Тогда $T$ продолжается на $\H_1$ 
как замкнутый ограниченный оператор из $\H_1$ в $\H_{-1}$. Из
условия~\eqref{III:eq:Qsub} и теоремы об устойчивости замкнутости 
(см.~\cite{Ka}) получаем, что $T+Q:\H_1\to\H_{-1}$ замкнут и ограничен,
кроме того, форма $((T+Q)x,x)$, определенная на $\H_1$, замкнута и 
секториальна. Согласно первой теореме о представлении, с этой формой
ассоциирован $m$--секториальный оператор $S=T\wt+Q$ в пространстве $\H$,
который есть сужение оператора $T+Q$ на область
$$
 \D(S)=\{x|Tx+Qx\in\H\}.
$$
Если операторы $Q_n$ таковы, что выплонено соотношение~\eqref{III:eq:Qconv},
то при достаточно больших $n$ для $Q_n$ выполняется 
неравенство~\eqref{III:eq:Qsub}, а потому определены операторы $T\wt+Q$. 
Докажем равномерную резольвентную сходимость~\eqref{III:eq:Sconv}.

Пусть $T_0$ --- оператор в $\H$, ассоциированный 
с квадратичной формой $\Re(Tx,x)$ (этот оператор уже был 
введен в начале параграфа). Так как $T_0$ равномерно положителен и 
$\D(T_0^{1/2})=\H_1$, то $T^{1/2}$ осуществляет изоморфизм $\H_1$ на 
$\H$ и $\H$ на $\H_{-1}$.

Поскольку нормы операторов $T+Q_n:\H_1\to\H_{-1}$ равномерно ограничены, 
то нормы операторов
\begin{equation}\label{III:eq:Zn}
Z_n=T_0^{-1/2}(T+Q_n)T_0^{-1/2}:\H\to\H
\end{equation}
также ограничены равномерно. Кроме того, из оценок~\eqref{III:eq:Qsub}
и~\eqref{III:eq:Qconv} следует, что 
$$
 c_1\|y\|^2_{\H}\le\Re((Z_n+\rho)y,y)\le c_2\|y\|^2_{\H}\quad\forall y\in\H
$$
с некоторыми константами $c_1$, $c_2>0$, не зависящими от $n$ и $y$,
если $\rho$ --- достаточно большое положительное число. Следовательно,
числовые образы операторов $Z_n+\rho$ отделены от нуля постоянной $c_1>0$.
Но тогда, согласно теореме Хаусдорфа, $Z_n+\rho$ обратимы, причем
$\|(Z_n+\rho)^{-1}\|\le c_1^{-1}$.

Теперь доказательство можно закончить так же, как в работе~\cite{MSh}. Имеем
\begin{multline*}
 (T+Q_n-\rho)^{-1}-(T+Q-\rho)^{-1}=(T+Q_n-\rho)^{-1}
(Q-Q_n)(T+Q-\rho)^{-1}=\\
=T_0^{-1/2}(Z_n+\rho)^{-1}
[T_0^{-1/2}(Q-Q_n)T_0^{-1/2}](Z+\rho)^{-1}T_0^{-1/2},
\end{multline*}
где $Z$ определен~\eqref{III:eq:Zn}, если в этом равенстве опустить индекс
$n$. Очевидно, последнее тождество доказывает равномерную резольвентную
сходимость~\eqref{III:eq:Sconv}. Теорема доказана.
\end{proof}
\begin{remark}\label{III:rm:spec}
В условиях теоремы~\ref{III:th:gensum} нельзя утверждать обычную
сходимость спектров операторов $T\wt+Q_n$ к спектру $T\wt+Q$. Контрпримеры
можно найти даже в случае ограниченных операторов $T$ и $Q$, действующих в
$\H$ (см., например,~\cite{Halmosh}). При каких дополнительных
условиях сходимость спектров утверждать можно? Укажем одно достаточное
условие: $T$ --- симметрический оператор, выполнена оценка~\eqref{III:eq:Qsub},
а $Q:\H_1\to\H_{-1}$ --- компактный оператор. В этом случае сходимость
спектров снизу (а значит, и обычную сходимость) можно доказать
точно так же, как и в самосопряженом случае.
\end{remark}

\subsection{Определение и свойства пространств $M[k,-l]$}

Пространство $H_2^s(\R^n)$ наделяется естественным скалярным 
произведением
$$
(f,g)_s=((1+x^2)^{s/2}Ff,(1+x^2)^{s/2}Fg)_{L_2}.
$$
Известно (см.~\cite{Tr}), что
оператор $(-\De+1)^{\al/2}$ осуществляет изоморфизм пространств $H_2^s$ и
$H_2^{s-\al}$.

Известно также~\cite{BIN}, что при целых $s\in\R^+$
эквивалентное скалярное произведение в $H_2^s$ можно задать выражением
$$
 (f,g)'_s=\suml_{|\al|\le s} (D^\al f,D^\al g)_{L_2}.
$$

Обозначение $(f,g)$ зарезервируем для обозначения действия
функционала $f\in H^{-s}_2$  на элемент $g\in H^{s}_2$
при всех $s\ge0$, а также действия $f\in D'$ на $g\in D$.

Далее в этом параграфе числа $k$, $l\ge0$ не предполагаются целыми.

Пусть $\phi\in D'$ --- обобщенная функция. Определим полуторалинейную
форму $(\phi f,g)$ на $D\times D$ как $(\phi,f\ov{g})$.
Назовем $\phi$ мультипликатором из $H_2^k$ в $H_2^{-l}$ (обозначая
$\phi\in M[k,-l]$), если форма $(\phi f,g)$ продолжается по непрерывности
на $H_2^k\times H_2^l$, иными словами, имеет место оценка
$$
(\phi f,g)\le C\|f\|_k \|g\|_l\qquad
\text{для всех $\phi$, $\psi\in D$}.
$$

В этом случае формой $(\phi f,g)$ определяется ограниченный оператор
$F:H^k_2\to H^{-l}_2$, а именно, каждому $f\in H^k_2$ ставим в
ссответствие функционал $M_\phi f\in H^{-l}_2$, действующий на $g\in
H^l_2$ как $(M_\phi f,g)=(\phi f,g)$.  Заметим, что $M_\phi$ является
оператором поточечного умножения на $\phi$, если $\phi\in D$.  Норма
элемента $\phi$ в пространстве мультипликаторов полагается равной норме
оператора $M_\phi$.

Выделим в $M[k,-l]$ подпространства $M_0[k,-l]$ и $\Mo[k,-l]$.
$\Mo[k,-l]$ определим как замыкание пространства основных функций $D$ по
норме $M[k,-l]$.

Подпространство $M_0[k,-l]$ по определению будет состоять из элементов
$\phi$ таких, что оператор $M_\phi$ является
$(-\De+1)^{\frac{k+l}2}$--ограниченным с произвольно малой относительной 
нормой,
иными словами, для произвольной $f\in H_2^k$ выполнена оценка
$$
\|M_\phi f\|_{-l}\le \ep
\|(-\De+1)^{\frac{k+l}2} f\|_{-l}+C(\ep)\|f\|_{k}.
$$

В дальнейшем будет использоваться интерполяционная теорема
для шкалы гильбертовых пространств $H_2^s$. Она нам понадобится в
следующей формулировке (см., нампример,~\cite{LM}):
\begin{theorem*}{Теорема об интерполяции}
Пусть задан линейный оператор $F:D\to D$, и имеют место оценки
\begin{align*}
 &\|F f\|_{s_1}\le C_1\|f\|_{t_1},\\
 &\|F f\|_{s_2}\le C_2\|f\|_{t_2}.
\end{align*}
Тогда при $0<\al<1$ выполнена оценка
$$
 \|F f\|_{\al s_1+(1-\al) s_2}\le C_1^{\al} C_2^{1-\al}\|f\|_{\al
 t_1+(1-\al) t_2}.
$$
\end{theorem*}

Перейдем к описанию свойств пространств $M[k,-l]$.

Сразу же заметим, что пространства $M[k,-l]$ и $M[l,-k]$ совпадают.
Действительно, оценки
$$
|(\phi,f\ov{g})|\le C\|f\|_k \|g\|_l\qquad\forall f,g\in D
$$
и
$$
|(\phi,f\ov{g})|\le C\|f\|_l \|g\|_k\qquad\forall f,g\in D
$$
выполняются одновременно.

Из вышесказанного и интерполяционной теоремы немедленно следует
\begin{lemm}
Если $k_2<k_1$ и $k_1+l_1=k_2+l_2$, то $M[k_1,-l_1]\subset M[k_2,-l_2]$,
и $\|\cdot\|_{M[k_2,-l_2]}\le \|\cdot\|_{M[k_2,-l_2]}$.
\end{lemm}

В следующем параграфе нам потребуется следующее утверждение.
\begin{lemm}
Если $\phi\in M[k,-l]$, то квадратичная форма $\phi_m[f]=
(\phi D^{m-k}f,D^{m-l}f)$ продолжается по непрерывности на
$H_2^m$. Если $\phi\in M_0[k,-l]$, то $\phi_m[f]$
подчинена форме $(\nabla^m f,\nabla^m f)$ с произвольно малой относительной
гранью.
\end{lemm}
\begin{proof}
Для $f$, $g\in D$ имеем оценку
$$
 (\phi D^{m-k}f,D^{m-l}g)\le \|\phi\|_{M[k,-l]} \|D^{m-k} f\|_k
 \|D^{m-l} g\|_l \le C \|f\|_m \|g\|_m,
$$
что доказывает первое утверждение. Пусть теперь $\phi\in
M_0[k,-l]$. Тогда
\begin{multline*}
|(\phi D^{m-k}f,D^{m-l}f)|\le \ep
((-\De+1)^{\frac{k+l}2}D^{m-k}f,D^{m-l}f)_s+C(\ep) (f,f)=\\
=\ep (Tf,f)+C(\ep) (f,f),
\end{multline*}
где $T=D^{m-l}(-\De+1)^{\frac{k+l}2}D^{m-k}$. Оператор $T$ ---
дифференциальный порядка $2m$, поэтому он ограничен как оператор из
$H^m_2$ в $H^{-m}_2$, т.е. $(Tf,f)\le C(f,f)_m$. Следовательно,
$$
|(\phi D^{m-k}f,D^{m-l}f)|\le \ep_1 (f,f)_m+C(\ep_1) (f,f).
$$
\end{proof}

Приведем ряд достаточных условий принадлежности функции пространству 
$M[k,-l]$.
Поскольку $M[k,-l]=M[l,-k]$, то
 далее без ограничения общности предположим $k>l$. В этом
предположении $M[k,-l]\subset H_{2,loc}^{-l}$.

\begin{lemm}\label{III:th:H2}
Если $k>n/2$, то $H_2^{-l}\subset \Mo[k,-l]$, причем
$\|\cdot\|_{M[k,-l]}\le C\|\cdot\|_{-l}$.
\end{lemm}
\begin{proof}
Для начала докажем следующее утверждение:

{\it
Пусть $k>n/2$. Тогда для любого $k\ge l\ge-k$ и произвольной пары функций
$f$, $g\in D$ выполнена оценка
$$
 \|fg\|_l\le C\|f\|_k \|g\|_l.
$$
}
Для $l=k$ это утверждение приведено в~\cite{MSh}. Для $l=-k$ оно доказано в
главе 1. Теперь для произвольного $l$ из $(-k,k)$ утверждение
получается путем применения интерполяционной теоремы к оператору
умножения на $f$.

Вернемся к доказательству предложения. Пусть $\phi\in H_2^{-l}$. В наших
предположениях для произвольных $f$, $g\in D$ верна оценка
$$
|(\phi f,g)|\le \|\phi\|_{-l} \|fg\|_l\le C \|\phi\|_{-l} \|f\|_k\|g\|_l,
$$
откуда $\|\phi\|_{M[k,-l]}\le C \|\phi\|_{-l}$. Поскольку $D$ плотно в
$H_2^{-l}$, выполнено включение $H^{-l}_2\subset \Mo[k,-l]$.
\end{proof}

\begin{lemm}\label{III:th:Hp}
Пусть $k\le n/2$, $\ga\le l$ и $p>\dfrac{n}{k+l-\ga}$.
Тогда $H^\ga_p\subset\Mo[k,-l]$
и $\|\cdot\|_{M[k,-l]}\le \|\cdot\|_{\ga,p}$.
\end{lemm}
\begin{proof}
Доказательство аналогично доказательству теоремы~5 статьи~\cite{BSh}. 
В качестве крайних мультипликативных оценок надо взять 
$$
 \|\phi\psi\|_{0,1}\le\|\phi\|_{0,2}\|\psi\|_{0,2}
$$
и
$$
 \|\phi\psi\|_{L,2}\le\|\phi\|_{K,2}\|\psi\|_{L,2},
$$
где $K=n/2+\ep$, $L=K\frac{l}{k}$. Первая оценка --- это неравенство Гельдера,
вторая доказана в лемме 3.3. Полилинейная интерполяционная теорема 
(см.~\cite{BL})
даст оценку
$$
 \|\phi\psi\|_{l,p'_1}\le\|\phi\|_{k,2}\|\psi\|_{l,2},
$$
где $p'_1=\dfrac{n+2\ep}{n+2\ep-k}$. Эта оценка доказывает вложение 
$H^{-l}_{p_1}\subset M[k,-l]$ при $p_1\ge \frac{n}{k}$. Вложение
$H^{-\ga}_p\subset M[k,-l]$ при $p\ge \dfrac{n}{l+k-\ga}$ следует из теоремы
вложения Соболева.
\end{proof}

\begin{lemm}\label{III:th:Polking}
Пусть $k<n/2$, $\ga\le l$ и $p\ge\dfrac{n}{k+l-\ga}$. 
Тогда $H^\ga_p\subset\Mo[k,-l]$
и $\|\cdot\|_{M[k,-l]}\le \|\cdot\|_{\ga,p}$.
\end{lemm}
\begin{proof}
Доказательство по сути аналогично доказательству теоремы~6 статьи~\cite{BSh}.
Перечислим отличия:

Проверяться будет оценка
$$
 \|\phi\psi\|_{l,s}\le\|\phi\|_{k,2}\|\psi\|_{l,2}
$$
для $s=n/(n-k)$, $k\le l$, $k<n/2$.
В качестве параметров $\eta_{1,2}$, $\rho_{1,2}$, $q_{1,2}$ берем
\begin{align*}
&0<\eta_1<\max\biggl\{\eta,\frac np-k\biggr\},
\qquad
\eta_2=\eta-\eta_1,
\\
&\rho_1=\frac {2(n-k)}{n+2(|j|+\eta_1-k)},
\quad
\rho_2=\frac {2(n-k)}{n-2(|j|+\eta_1)},
\\
&q_1=2\rho_1,
\quad
q_2=2\rho_2.
\end{align*}
Неравенство
$$
\|\phi\|_{|j|+\eta_1,\rho_1 s}\|\psi\|_{|m|+\eta_2,\rho_2 s}
\le C\|\phi\|_{k,2}\|\psi\|_{l,2}
$$
следует из теоремы Соболева ввиду неравенств
$$
 k-(|j|+\eta_1)\ge \frac n{\rho_1 s}-\frac n2,\quad
 l-(|m|+\eta_2)\ge \frac n{\rho_2 s}-\frac n2,
$$
проверяемых непосредственно.
\end{proof}

\begin{lemm}\label{III:th:Fubini}
Пусть $\phi$ зависит от $p$ переменных и принадлежит пространству 
$M[k,-l]$ в $\R^p$, а функция $\psi$ зависит от $q$ переменных и
принадлежит $L_\infty(\R^q)$. Тогда функция $\phi(x)\psi(y)\in M[k,-l]$ в 
$\R^{p+q}$ и справедлива оценка
$$
 \|\phi\psi\|_{M[k,-l]}\le \|\phi\|_{M[k,-l]}\|\psi\|_{L_\infty}.
$$
Если же $\phi\in M_0[k,-l]$, то и $\phi\psi\in M_0[k,-l]$ в $\R^{p+q}$.
\end{lemm}
\begin{proof}
Дословно повторяет доказательство теоремы~8 статьи~\cite{NeSh}.
\end{proof}

\begin{theorem}\label{III:th:main}
Пусть оператор $L$ вида~\eqref{III:eq:L} таков, что для его главной части
$L_0$ выполнена коэрцитивная оценка~\eqref{III:eq:Re>}
(например, $L_0$ сильно эллиптичен и его коэффициенты непрерывны).
Пусть функции
$c_{\al,\be}$ при $|\al|+|\be|<2m$ принадлежат пространствам
$M_0[m-|\al|,-(m-|\be|)]$. Тогда оператор
$L$ корректно определен как
секториальный оператор в $L_2(\R^n)$, причем замыкание области определения
квадратичной формы этого оператора совпадает с $H_2^m$.

$H^m_2$. Если
последовательности гладких функций $c_{\al,\be,n}\in D$, $|\al|+|\be|<2m$,
сходятся к
$c_{\al,\be}$ по норме $M[m-|\al|,-(m-|\be|)]$, то операторы $L_n$ с
коэффициентами $c_{\al,\be,n}$
сходятся к оператору $L$ в смысле равномерной резольвентной сходимости.
При этом спектры $\si(L_n)$ сходятся к спектру $\si(L)$ сверху. Если же
$L_0$ --- симметрический оператор, а $c_{\al,\be}$ --- компактные
мультипликаторы в указанных пространствах, то имеет место обычная сходимость
спектров.
\end{theorem}
\begin{proof}
Эта теорема следует из теоремы~\ref{III:th:gensum} и
замечания~\ref{III:rm:spec}.
\end{proof}

Следующее утверждение можно считать основным результатом этой статьи. 
В формулировке используется понятие сходимости спектров сверху и снизу.
Соответствующие определения можно найти в книге~\cite{Ka}.

\begin{theorem}
Пусть главная часть $L_0$ оператора $L$ есть равномерно сильно
эллиптический оператор с непрерывными коэффициентами, а коэффициенты
$c_{\al,\be}$ являются распределениями, причем
\begin{align*}
 &c_{\al,\be}\in H^{|\al|-m}_p\quad
\text{при $|\al|\ge|\be|$, $p>\max\{2,\dfrac{n}{m-|\be|}\}$},\\
 &c_{\al,\be}\in H^{|\be|-m}_p\quad
\text{при $|\al|\le|\be|$, $p>\max\{2,\dfrac{n}{m-|\al|}\}$}.
\end{align*}
В случае $m-|\al|\ne n/2$ ($m-|\be|\ne n/2$) допускаются также значения
$p=\max\{2,\dfrac{n}{m-|\be|}\}$ ($p=\max\{2,\dfrac{n}{m-|\al|}\}$).
Тогда оператор $L$ корректно определен как $m$--секториальный оператор в
$\H$. Если $c_{\al,\be}^n(x)$ --- гладкие функции, сходящиеся в указанных 
пространствах $H^{|\al|-m}_p$ к $c_{\al,\be}$, то $L$ является равномерным
резольвентным пределом соответствующих регулярных операторов $L_n$. Спектры
$\si(L_n)$ сходятся к $\si(L)$ сверху, а если главный оператор $L_0$ 
симметричен, то имеет место и сходимость спектров снизу.
\end{theorem}
\begin{proof}
Это утверждение есть следствие леммы~\ref{III:th:Hp} и предыдущей теоремы.
Нужно лишь заметить, что умножение на гладкую функцию с компактными
носителем есть компактный мультипликатор. Но тест--функции плотны в
соболевских пространствах $H_p^{|\al|-m}$, поэтому умножение на функцию 
из этих классов будет также компактным мультипликатором в 
$M[m-|\al|,-(m-|\be|)]$, если число $p$ удовлетворяет условиям теоремы. 
Теорема доказана.
\end{proof}

\end{document}